\documentclass[12pt]{amsart}

\usepackage[all]{xy}

\theoremstyle{plain}
\newtheorem{theorem}{Theorem}[section]
\newtheorem{corollary}[theorem]{Corollary}
\newtheorem{lemma}[theorem]{Lemma}
\newtheorem{proposition}[theorem]{Proposition}
\newtheorem{question}[theorem]{Question}

\theoremstyle{definition}
\newtheorem{definition}[theorem]{Definition}
\newtheorem{example}[theorem]{Example}

\newtheorem{remark}[theorem]{Remark}
\newtheorem{Definitions and Notation}[theorem]{Definitions and
Notation}

\numberwithin{equation}{section}

\newcommand{\f}[1]{\ensuremath{\mathfrak{#1}}}
\newcommand{\locoho}[3]{\ensuremath{\operatorname{H}_{#1}^{#2}\left(#3\right)}}
\newcommand{\intersect}{\cap}
\newcommand{\lra}{\longrightarrow}
\newcommand{\sse}{\subseteq}

\newcommand{\ol}[1]{\ensuremath{\overline{#1}}}

\newcommand{\soc}[1]{\ensuremath{\operatorname{Soc}\left(#1\right)}}
\newcommand{\socdim}[1]{\ensuremath{\operatorname{Socdim}\left(#1\right)}}
\newcommand{\koszulh}[3]{\ensuremath{\operatorname{H}^{#1}\left(#2; #3\right)}}
\newcommand{\koszul}[3]{\ensuremath{\operatorname{K}^{#1}\left(#2; #3\right)}}

\renewcommand{\index}[3]{\ensuremath{\operatorname{N}_{#1}\left(#2;#3\right)}}
\newcommand{\seq}[3]{\ensuremath{{#1}_{#2}}, \ldots, \ensuremath{{#1}_{#3}}{}}
\newcommand{\undseq}[1]{\ensuremath{\underline{#1}}}
\newcommand{\directlimit}[2]{\ensuremath{\varinjlim_{#1} #2}}
\newcommand{\ideal}[4]{\ensuremath{\left({#1}_{#2}, \ldots, {#1}_{#3}\right)#4}}
\newcommand{\leaveout}[1]{\ensuremath{\widehat{#1}}}
\renewcommand{\hom}[3]{\ensuremath{\operatorname{Hom}_{#1}\left(#2,#3\right)}}
\renewcommand{\colon}[3]{\ensuremath{\left(#1:_{#2} #3\right)}}
\newcommand{\ext}[4]{\ensuremath{\operatorname{Ext}_{#1}^{#2}\left(#3,#4\right)}}
\newcommand{\tor}[4]{\ensuremath{\operatorname{Tor}_{#1}^{#2}\left(#3,#4\right)}}

\newcommand{\unmixed}[1]{\ensuremath{\operatorname{U}\left(#1 \right)}}

\newcommand{\depth}{\operatorname{depth}}
\newcommand{\Dim}{\operatorname{dim}}
\newcommand{\Ass}{\operatorname{Ass}}

\begin{document}

\title[Index of Reducibility]{The index of reducibility of parameter ideals and mostly zero finite
local cohomologies}

\author{Jung-Chen Liu}
\address{Department of Mathematics, National Taiwan Normal University, Taipei, Taiwan}
\email{liujc@math.ntnu.edu.tw}

\author{Mark W. Rogers}
\address{Department of Mathematics, Southwest Missouri State University, Springfield, MO 65805}
\email{markrogers@smsu.edu}

\begin{abstract}
In this paper we prove that if $M$ is a finitely-generated module of
dimension $d$ with finite local cohomologies over a Noetherian local
ring $(A, \f{m})$, and if $\locoho{\f{m}}{i}{M} = 0$ except possibly
for $i \in \{ 0, r, d \}$ with some $0 \leq r \leq d$, then there
exists an integer $\ell$ such that every parameter ideal for $M$
contained in $\f{m}^{\ell}$ has the same index of reducibility. This
theorem generalizes earlier work of the second author and is closely
related to recent work of Goto-Suzuki and Goto-Sakurai; Goto-Sakurai
have supplied an answer of \emph{yes} in case $M$ is Buchsbaum.
\end{abstract}



\maketitle

\section{Introduction}

Throughout this paper, $A$ denotes a Noetherian local ring with
maximal ideal \f{m} and residue field $k = A/\f{m}$, and $M$ is a
finitely generated $A$-module of dimension $d$.  In the same 1921
paper where she introduced the fundamental notion of primary
decomposition of ideals, E. Noether showed that every ideal in a
Noetherian ring may be expressed as a finite intersection of
irreducible ideals.  She also showed that if the intersection is
irredundant, then the number of irreducible ideals appearing in the
intersection depends only on the ideal and not on the particular
intersection.  This result readily generalizes to Noetherian
modules.  For a submodule $N$ of $M$, we refer to the number of
irreducible submodules appearing in any irredundant expression of
$N$ as an intersection of irreducible submodules as the \emph{index
of reducibility} of $N$.

By a \emph{parameter ideal} for $M$ we mean an ideal \f{q} that can
be generated by $d$ elements and such that $M/\f{q}M$ has finite
length.  We define the \emph{index of reducibility of \f{q}} on $M$
to be the index of reducibility of the submodule $\f{q}M$; we denote
the index of reducibility of \f{q} on $M$ by \index{A}{\f{q}}{M}.
Since $M/\f{q}M$ has finite length, it is known that the index of
reducibility of a parameter ideal \f{q} on $M$ is given by the
\emph{socle dimension} of $M/\f{q}M$; that is, $\index{A}{\f{q}}{M}
= \Dim_k \hom{A}{k}{M}$.  (The socle of $M$, which we denote by
\soc{M}, is generally defined as the sum of all the simple
submodules of $M$; in our setting, it equals \colon{0}{M}{\f{m}} and
is isomorphic to \hom{A}{k}{M}.)

In 1956, D. G. Northcott showed that the index of reducibility of
a parameter ideal in a Cohen-Macaulay local ring depends only on
the ring and not on the parameter ideal \cite[Theorem~3]{Northc}.
This result extends to modules, and the common index of
reducibility of parameter ideals on a Cohen-Macaulay module is
called the \emph{type} of the module. Regardless of whether $M$ is
Cohen-Macaulay, the type is defined to be $\Dim_{k}
\ext{A}{t}{k}{M}$, where $t$ is the depth of $M$.

Although Northcott and D.G. Rees proved in 1956 that if every
parameter ideal of a Noetherian local ring is irreducible then the
ring is Cohen-Macaulay \cite{NR-princ}, the constant index of
reducibility of parameter ideals does not characterize
Cohen-Macaulay local rings. Indeed, in 1964 S. Endo and M. Narita
gave examples of non-Cohen-Macaulay Noetherian local rings having
constant index of reducibility of parameter ideals
\cite{Endo.Narita}.

In 1984, S. Goto and N. Suzuki revived interest in the index of reducibility of
parameter ideals, generalizing the examples of Endo-Narita, and paying
particular attention to the supremum of the index of reducibility of parameter
ideals \cite{Goto.Suzuki}. We refer to this supremum as the \emph{Goto-Suzuki
type}. Although examples were given showing that the Goto-Suzuki type may be
infinity, finite upper and lower bounds were provided in the case where $M$ has
finite local cohomologies \cite[Theorem~2.1 and Theorem~2.3]{Goto.Suzuki}.

\begin{definition}
We say that $M$ has \emph{finite local cohomologies} if the local cohomology $A$-modules
\locoho{\f{m}}{i}{M} have finite length for each integer $i$ such that $0 \leq i \leq d - 1$.
\end{definition}

Implicit in \cite{Goto.Suzuki} we find that if $M$ is a module
having finite local cohomologies, there exists inside every power of
\f{m} a parameter ideal \f{q} such that
\[
\index{A}{\f{q}}{M} = \sum_{i = 0}^{d} \binom{d}{i} \socdim{\locoho{\f{m}}{i}{M}}.
\]
This is the lower bound for the Goto-Suzuki type mentioned above.

An equivalent condition for $M$ to have finite local cohomologies
is the existence of a standard ideal for $M$.
\begin{definition}
An ideal \f{a} is called a \emph{standard ideal} for $M$ if \f{a}
is an \f{m}-primary ideal with the property that for each system
of parameters $x_1$, \ldots, $x_d$ of $M$ contained in \f{a}, we
have
\[
\colon{\ideal{x}{1}{i - 1}{M}}{M}{x_{i}} = \colon{\ideal{x}{1}{i -
1}{M}}{M}{\f{a}}
\]
for every integer $i$ with $1 \leq i \leq d$.
\end{definition}

In 2003, Goto and H. Sakurai showed that if $M$ is Buchsbaum (i.e.,
\f{m} is a standard ideal for $M$), then there exists a power of
\f{m} inside which every parameter ideal for $M$ has the same index
of reducibility on $M$ \cite[Corollary~3.13]{Goto.SakuraiI},
necessarily equal to the lower bound of the Goto-Suzuki type.  We
refer to this by saying $M$ has \emph{eventual constant index of
reducibility of parameter ideals}.  In a recent paper, the second
author showed that if $M$ has dimension one, or if $M$ has dimension
two, finite local cohomologies, and positive depth, then $M$ has
eventual constant index of reducibility of parameter ideals
\cite[Theorem~2.3 and Theorem~3.3]{Rogers}.  The results mentioned
in this paragraph provide partial answers to the following question,
which appears as \cite[Question~1.2]{Rogers}.

\begin{question}
Suppose $(A,~\f{m},~k)$ is a Noetherian local ring having finite
local cohomologies.  Does $A$ have eventual constant index of
reducibility of parameter ideals?
\end{question}

In this paper we provide a partial answer to this question by
generalizing the results in \cite{Rogers}. Our main result is the
following theorem.

\begin{theorem}[Main Theorem]\label{MainTheorem}\hfill\\
Suppose $M$ has finite local cohomologies and that
$\locoho{\f{m}}{i}{M} = 0$ except possibly for $i \in \{ 0, r, d
\}$, where $r$ is some integer with $0 \leq r \leq d$. There exists
an integer $\ell$ such that for every parameter ideal \f{q} of $M$,
if $\f{q} \sse \f{m}^{\ell}$ then the index of reducibility of \f{q}
on $M$ is independent of \f{q} and is given by
\[
\index{A}{\f{q}}{M} = \sum_{i = 0}^{d} \binom{d}{i}
\socdim{\locoho{\f{m}}{i}{M}}
\]
\end{theorem}
In the process of proving this result, we provide information on
the behavior of the unmixed components of parts of systems of
parameters contained in a standard ideal.  As a corollary to the
Main Theorem, we present a result that characterizes a certain
class of finite local cohomology modules as Gorenstein precisely
when every power of the maximal ideal contains an irreducible
parameter ideal (here we take as definition that a Gorenstein
local ring is a Cohen-Macaulay local ring having an irreducible
parameter ideal.)

\section{Proof of the main theorem and a corollary}

The importance of the existence of a standard ideal is made clear
by the following proposition, for which we refer the reader to
\cite[Corollary~18, p.~264]{Stuckrad.Vogel}.

\begin{proposition}\label{P:FLC=standardideal}
There is a standard ideal for $M$ if and only if $M$ has finite
local cohomologies.
\end{proposition}

Now we present the proof of the main theorem,
Theorem~\ref{MainTheorem}; the statement and proof of several
useful lemmas are postponed.

\begin{proof}
Set $W = \locoho{\f{m}}{0}{M}$ and $\ol{M} = M/W$.  According to
Proposition~\ref{P:positivedepth} there is an integer $n$ so that
for every parameter ideal \f{q} for $M$, if $\f{q} \sse \f{m}^{n}$
then $\index{A}{\f{q}}{M} = \socdim{M} + \index{A}{\f{q}}{\ol{M}}$.
Since $\locoho{\f{m}}{i}{\ol{M}} \cong \locoho{\f{m}}{i}{M}$ for all
$i > 0$, we may require the integer $\ell$ to be at least as large
as $n$ and replace $M$ by \ol{M} for the remainder of the proof to
assume that $\depth M > 0$. If $\locoho{\f{m}}{i}{M} = 0$ for all $i
< d$, then we're done since $M$ is Cohen-Macaulay. Otherwise, we
have that $\locoho{\f{m}}{i}{M} = 0$ for all $i$ except $i = d$ and
$i = r$ where $0 < r < d$.

For a submodule $N$ of $M$, we denote by \unmixed{N} the unmixed component of $N$; that is,
\unmixed{N} is the intersection of the primary components of $N$ whose associated primes have
maximal dimension, equal to $\Dim M/N$.  According to \cite[Proposition~3.2]{Rogers}, if we let
\f{a} be a standard ideal for $M$ then there exists an integer $\ell$ such that $\f{m}^{\ell} \sse
\f{a}$ and such that for any parameter ideal $\f{q} = \ideal{x}{1}{d}{A}$ for $M$, if $\f{q} \sse
\f{m}^{\ell}$ then the index of reducibility of \f{q} on $M$ is given by
\[
\index{A}{\f{q}}{M} = \socdim{
    \sum_{i = 1}^{d}
        \frac{ U_{i} + \f{q} M }{ \f{q}M }
    }
    + \socdim{\locoho{\f{m}}{d}{M}}.
\]
where $U_{i} = \unmixed{ \left( x_{1}, \ldots, \leaveout{x_{i}}, \ldots, x_{d} \right) M }$.  Our
task in this proof is to examine the sum appearing in the expression above.

We begin by using the assumption that $\locoho{\f{m}}{i}{M} = 0$ for all
integers $i$ with $r < i < d$ to apply Proposition~\ref{P:BustUpUnmixed} and
obtain
\[
\sum_{i = 1}^{d}
    \frac{U_{i} + \f{q}M}{\f{q}M}
= \sum_{1 \leq i_{1} < \cdots < i_{r} \leq d}
\frac{\unmixed{\left( x_{i_{1}}, \ldots, x_{i_{r}} \right) M} +
\f{q}M}{\f{q}M}.
\]

Next we examine one of the summands on the right side. Using
Lemma~\ref{P:UnmixedIntersectParameterIdeal} we see that
\[
\frac{
    \unmixed{
        \left(
            x_{i_{1}}, \ldots, x_{i_{r}}
        \right) M
    }
    +
    \f{q}M
}{
    \f{q}M
} \cong
\frac{\unmixed{\left( x_{i_{1}}, \ldots,
x_{i_{r}} \right) M}}{\left( x_{i_{1}}, \ldots, x_{i_{r}} \right)
M}.
\]
Since $\locoho{\f{m}}{i}{M} = 0$ for each integer $i$ such that $0 \leq i \leq r - 1$, we use
Proposition~\ref{P:Locoho=UnmixedModParameter} to obtain
\[
\frac{
    \unmixed{
        \left(
            x_{i_{1}}, \ldots, x_{i_{r}}
        \right) M
    }
    +
    \f{q}M
}{
    \f{q}M
} \cong \locoho{\f{m}}{r}{M}.
\]
It remains to show that the sum
\[
\sum_{1 \leq i_{1} < \cdots < i_{r} \leq d} \frac{\unmixed{\left(
x_{i_{1}}, \ldots, x_{i_{r}} \right) M} + \f{q}M}{\f{q}M}
\]
is direct.

By symmetry, it suffices to show that $\f{q} M$ equals the module
\[
\left[ \unmixed{\ideal{x}{1}{r}{M}} + \f{q}M \right] \intersect
\left[ \sum_{\substack{1 \leq i_{1} < \cdots < i_{r} \leq d\\
i_{r} > r}} \unmixed{ \left( x_{i_{1}}, \ldots, x_{i_{r}} \right) M}
+ \f{q}M \right].
\]
Using basic properties of intersections, this is equivalent to
showing that
\[
\unmixed{\ideal{x}{1}{r}{M}} \intersect \left[ \sum_{\substack{1 \leq i_{1} < \cdots < i_{r} \leq d\\
i_{r} > r}} \unmixed{ \left( x_{i_{1}}, \ldots, x_{i_{r}} \right)
M}  \right] \sse \f{q}M.
\]
Using Lemma~\ref{L:UnmixedEqualsColon}, we see that each summand
\unmixed{ \left( x_{i_{1}}, \ldots, x_{i_{r}} \right) M}, where
$i_{r} > r$, is contained in at least one of the submodules
\[
\unmixed{ \left( x_{1}, \ldots, \leaveout{x_{i}}, \ldots, x_{d}
\right) M}
\]
$(1 \leq i \leq r)$.  Thus the intersection above is contained in
\[
\unmixed{\ideal{x}{1}{r}{M}} \intersect \left[ \sum_{i = 1}^{r}
\unmixed{ \left( x_{1}, \ldots, \leaveout{x_{i}}, \ldots, x_{d}
\right) M} \right].
\]
According to Proposition~\ref{P:SumIsDirect}, this intersection is
contained in the submodule \ideal{x}{1}{r}{M}.  This completes the
proof of the theorem.
\end{proof}

\begin{corollary}
Suppose $A$ is a Noetherian local ring having finite local
cohomologies and that \locoho{\f{m}}{i}{A} is zero except possibly
for $i \in \{ 0, r, d \}$, where $r$ is some integer with $0 \leq r
\leq d$, and $d = \Dim A$. Then $A$ is Gorenstein if and only if
every power of \f{m} contains an irreducible parameter ideal.
\end{corollary}

\begin{proof}
It is well known that if $A$ is Gorenstein then every parameter
ideal is irreducible.  Now suppose that every power of \f{m}
contains an irreducible parameter ideal.  According to the Main
Theorem, each parameter ideal inside a high enough power of the
maximal ideal has index of reducibility
\[
\sum_{i = 0}^{d} \binom{d}{i} \socdim{\locoho{\f{m}}{i}{A}};
\]
this expression must equal 1. Since \locoho{\f{m}}{d}{A} is a
nonzero Artinian module, it has a nonzero socle.  Thus the socles of
\locoho{\f{m}}{i}{A} ($i < d$) must be zero, so that the Artinian
modules \locoho{\f{m}}{i}{A} are themselves zero. Since $A$ has only
one nonzero local cohomology module, $A$ is Cohen-Macaulay. Since
$A$ is Cohen-Macaulay with an irreducible parameter ideal, $A$ is
Gorenstein.  This completes the proof of the corollary.
\end{proof}

\section{Supporting propositions}

This section contains the supporting propositions
used in the proof of the Main Theorem.

\subsection{Reduction to positive depth}

\begin{proposition}\label{P:positivedepth}
Let $W = \locoho{\f{m}}{0}{M}$.  There exists an integer $n$ such
that for every parameter ideal \f{q} for $M$, if $\f{q} \sse
\f{m}^n$ then the index of reducibility of \f{q} on $M$ is given by
\[
\index{A}{\f{q}}{M} = \socdim{M} + \index{A}{\f{q}}{M/W}.
\]
\end{proposition}

\begin{proof}
Since $W$ has finite length, there is a positive integer $a$ such
that $W \cap \f{m}^{a} M = 0$.  Set $\ol{M} = M/W$.  Since \ol{M}
has positive depth, there is an \ol{M}-regular element $x$ of $A$.
By the Artin-Rees Lemma, there is a positive integer $b$ such that
for all positive integers $c$, we have $\f{m}^{b+c} \ol{M} \cap x
\ol{M} = \f{m}^{c} \left( \f{m}^{b} \ol{M} \cap x \ol{M} \right)$.
Let $n = a + b$ and suppose \f{q} is a parameter ideal for $M$ that
is contained in $\f{m}^{n}$.

Since $\f{q} M \cap W = 0$, we have $(W + \f{q} M)/\f{q} M \cong W$,
and thus an exact sequence
\[
0 \lra W \lra M / \f{q} M \lra \ol{M} / \f{q} \ol{M} \lra 0.
\]
Since $\soc{W} = \soc{M}$, after applying the socle functor $\soc{*}
= \colon{0}{*}{\f{m}}$ we obtain the exact sequence
\[
0 \lra \soc{M} \lra \soc{M / \f{q} M} \lra \soc{\ol{M} / \f{q}
\ol{M}}.
\]
Using the additivity of length, we see that the proof will be
complete if we show that the map $\soc{M / \f{q} M} \lra \soc{\ol{M}
/ \f{q} \ol{M}}$ is surjective.

Let $s \in M$ be a representative of a nonzero element of
\[
\soc{\ol{M} / \f{q} \ol{M}} \cong \colon{(\f{q} M + W)}{M}{\f{m}} /
(\f{q} M + W).
\]
Let \ol{s} denote the class of $s$ in \ol{M}. Then $\ol{s} \in
\colon{\f{q} \ol{M}}{\ol{M}}{\f{m}}$. Note that since
\[
x \colon{\f{m}^{n} \ol{M}}{\ol{M}}{\f{m}} \sse x \colon{\f{m}^{n}
\ol{M}}{\ol{M}}{x} = \f{m}^{n} \ol{M} \cap x \ol{M} = \f{m}^{a}
\left( \f{m}^{b} \ol{M} \cap x \ol{M} \right) \ol{M},
\]
so that $x \colon{\f{m}^{n} \ol{M}}{\ol{M}}{\f{m}} \sse x
\f{m}^{a}$.  Since $x$ is \ol{M}-regular, we have
\[
\colon{\f{m}^{n} \ol{M}}{\ol{M}}{\f{m}} \sse \f{m}^{a} \ol{M}.
\]
Since $\f{q} \sse \f{m}^{n}$, we see that $\ol{s} \in \f{m}^{a}
\ol{M}$, so that $s \in \f{m}^{a} M + W$.  Replace $s$ by another
representative of \ol{s} so that we may assume $s \in \f{m}^{a} M$.

Now we have $\f{m} s \sse \f{m}^{a + 1} M \cap (\f{q} M + W)$. Since
$\f{q} M \sse \f{m}^{a + 1} M$ and $\f{m}^{a + 1} M \cap W = 0$, we
have $\f{m} s \sse \f{q} M$.  Thus $s \in \soc{M / \f{q} M}$.  This
shows that the map $\soc{M / \f{q} M} \lra \soc{\ol{M} / \f{q}
\ol{M}}$ is surjective, so our proof is complete.
\end{proof}

\subsection{The unmixed component of a sum versus the sum of the unmixed components}

\begin{definition} Let $N$ be a submodule of $M$.
We say that $N$ is \emph{unmixed up to \f{m}-primary components}
if $\f{p} \in \Ass M/N$ implies that either $\f{p} = \f{m}$ or
$\Dim A/\f{p} = \Dim M/N$.
\end{definition}

\begin{remark}\label{R:Unmixed}\hfill
\begin{enumerate}
\item The unmixed component does not depend on a particular
primary decomposition.

\item If $N$ is a submodule of $M$ that is unmixed up to
\f{m}-primary components, then $\locoho{\f{m}}{0}{M/N} \cong
\unmixed{N}/N$.

\item According to \cite[Lemma~2.2, p.~71]{Stuckrad.Vogel}, if $M$ has finite local
cohomologies then for any part $x_1$, \ldots, $x_r$ $(r \geq 0)$
of a system of parameters for $M$, the submodule
\ideal{x}{1}{r}{M} is unmixed up to \f{m}-primary components.
\end{enumerate}
\end{remark}

The following lemma explains how we usually think of the unmixed
component of a submodule generated by part of a system of
parameters contained in a standard ideal.

\begin{lemma}\label{L:UnmixedEqualsColon}
Suppose $M$ has finite local cohomologies and let \f{a} be a
standard ideal for $M$.  Suppose \seq{x}{1}{r} $(1 \leq r \leq d)$
is part of a system of parameters for $M$.  If $\ideal{x}{1}{r}{A}
\sse \f{a}$ then
\[
\unmixed{\ideal{x}{1}{r - 1}{M}} = \colon{\ideal{x}{1}{r -
1}{M}}{M}{x_r}.
\]
\end{lemma}

\begin{proof}
Since \ideal{x}{1}{r - 1}{M} is unmixed up to \f{m}-primary components, there is an integer $n$
such that
\[
\unmixed{\ideal{x}{1}{r - 1}{M}} = \colon{\ideal{x}{1}{r - 1}{M}}{M}{x_{r}^{n}}.
\]
This last expression equals \colon{\ideal{x}{1}{r - 1}{M}}{M}{x_r}, since they are both equal to
\colon{\ideal{x}{1}{r - 1}{M}}{M}{\f{a}}.
\end{proof}

The following lemma is essentially a collection of several results
from \cite{Stuckrad.Vogel}.

\begin{lemma}\label{L:ColonToUnmixeds}
Suppose $M$ has finite local cohomologies and let \f{a} be a
standard ideal with respect to $M$.  Suppose \seq{x}{1}{r} $(1
\leq r \leq d)$ is part of a system of parameters for $M$.  If
$\ideal{x}{1}{r}{A} \sse \f{a}$, then for all integers $n_1$,
\ldots, $n_r \geq 1$ we have
\begin{align*}
& \colon{(x_{1}^{n_{1} + 1}, \ldots, x_{r}^{n_{r} +
1})M}{M}{x_{1}^{n_{1}} \cdot \cdots
\cdot x_{r}^{n_{r}}} \\
& \qquad = \ideal{x}{1}{r}{M} + \sum_{i = 1}^{r} \unmixed{(x_1,
\ldots, \widehat{x_i}, \ldots, x_{r})M}.\notag
\end{align*}
\end{lemma}

\begin{proof}
According to \cite[Theorem~20, Lemma~23, Lemma~24, pp
264--266]{Stuckrad.Vogel}, we have
\begin{align*}
& \colon{(x_{1}^{n_{1} + 1}, \ldots, x_{r}^{n_{r} +
1})M}{M}{x_{1}^{n_{1}} \cdot \cdots
\cdot x_{r}^{n_{r}}} \\
& \qquad = \ideal{x}{1}{r}{M} + \sum_{i = 1}^{r} \colon{(x_1,
\ldots, \widehat{x_i}, \ldots, x_{r})M}{M}{x_i}.\notag
\end{align*}
We complete the proof with an application of
Lemma~\ref{L:UnmixedEqualsColon}.
\end{proof}

Now we recall the connection between Koszul cohomology and local
cohomology.  Let $I$ be an ideal generated by elements
\seq{y}{1}{r}.  For a positive integer $n$, we use $\undseq{y}^n$
to denote the sequence $y_{1}^{n}$, \ldots, $y_{r}^{n}$.  We use
\koszulh{i}{\undseq{y}}{M} to denote the $i$th cohomology module
of the Koszul cocomplex \koszul{\bullet}{\undseq{y}}{M}.  Since
there are containments of ideals
\[
(\undseq{y})A \supseteq (\undseq{y}^2)A \supseteq (\undseq{y}^3)A
\supseteq \cdots,
\]
for each $i$ the corresponding Koszul cohomology modules fit into
a direct system
\[
\koszulh{i}{\undseq{y}}{M} \lra \koszulh{i}{\undseq{y}^2}{M} \lra
\koszulh{i}{\undseq{y}^3}{M} \lra \cdots.
\]
When $i = r$, the differentials are particularly simple:  they are
each multiplication by the product $y_1 \cdot \cdots \cdot y_r$.

For each $i$, it is known that the direct limit of the Koszul
cohomology is isomorphic to the cohomology module
\locoho{I}{i}{M}; that is,
\[
\directlimit{n}{\koszulh{i}{\undseq{y}^n}{M}} \cong
\locoho{I}{i}{M}.
\]

\begin{lemma}\label{L:UnmixedInsideColon}
Suppose $M$ has finite local cohomologies and
$\locoho{\f{m}}{r}{M} = 0$ for some $1 \leq r \leq d - 1$.  If
$x_1$, \ldots, $x_{r+1}$ is part of a system of parameters for
$M$, then there exists an integer $n$ such that
\[
\colon{(x_{1}, \ldots, x_{r})M}{M}{x_{r+1}} \sse
\colon{(x_{1}^{n+1}, \ldots, x_{r}^{n+1})M}{M}{(x_1 \cdots
x_r)^n}.
\]
\end{lemma}

\begin{proof}
Let $\f{q}_r = \ideal{x}{1}{r}{A}$ and $\f{q}_{r+1} =
\ideal{x}{1}{r + 1}{A}$. According to \cite[Lemma~1.5,
p.~28]{Stuckrad.Vogel}, we have the following commutative diagram:
\[
\xymatrix{
 \koszulh{r}{\seq{x}{1}{r + 1}}{M} \ar[r] \ar[d] &
\koszulh{r}{\seq{x}{1}{r}}{M} \ar[d]\\
 \locoho{\f{q}_{r + 1}}{r}{M} \ar[r] & \locoho{\f{q}_r}{r}{M}
 }
\]
According to \cite[Lemma~22, p.~264]{Stuckrad.Vogel}, we have
$\locoho{\f{q}_{r+1}}{r}{M} \cong \locoho{\f{m}}{r}{M}$, and by
hypothesis, this module is zero. (The application of this lemma is
the only place in the proof that we need finite local
cohomologies.) Thus, our diagram becomes:
\[
\xymatrix{
 \koszulh{r}{\seq{x}{1}{r + 1}}{M} \ar[r] \ar[d] &
\koszulh{r}{\seq{x}{1}{r}}{M} \ar[d]\\
 0 \ar[r] & \locoho{\f{q}_r}{r}{M}
 }
\]
Now, we apply the covariant functor $\koszulh{0}{x_{r + 1}}{-}$,
which simply returns the annihilator of $x_{r + 1}$:
\[
\xymatrix{
 \koszulh{0}{x_{r + 1}}{\koszulh{r}{\seq{x}{1}{r + 1}}{M}}
\ar[r] \ar[d] & \koszulh{0}{x_{r +
1}}{\koszulh{r}{\seq{x}{1}{r}}{M}} \ar[d]\\
 0 \ar[r] & \koszulh{0}{x_{r + 1}}{\locoho{\f{q}_r}{r}{M}}
  }
\]
Since $x_{r + 1}$ annihilates \koszulh{r}{\seq{x}{1}{r + 1}}{M}, we have
\[
\koszulh{0}{x_{r + 1}}{\koszulh{r}{\seq{x}{1}{r + 1}}{M}} = \koszulh{r}{\seq{x}{1}{r + 1}}{M}.
\]
From \cite[Corollary~1.7, p.~29]{Stuckrad.Vogel} we know that the top row in our diagram has become
surjective. Thus we have the following commutative diagram:
\[
\xymatrix{
 \koszulh{r}{\seq{x}{1}{r + 1}}{M} \ar[r] \ar[d] &
\koszulh{0}{x_{r + 1}}{\koszulh{r}{\seq{x}{1}{r}}{M}} \ar[d]
\ar[r] & 0\\
 0 \ar[r] & \koszulh{0}{x_{r + 1}}{\locoho{\f{q}_r}{r}{M}}
  }
\]
From a simple diagram chase, we see that the map
\[
\koszulh{0}{x_{r + 1}}{\koszulh{r}{\seq{x}{1}{r}}{M}} \lra
\koszulh{0}{x_{r + 1}}{\locoho{\f{q}_r}{r}{M}}
\]
is the zero map. Thus the submodule
\[
\koszulh{0}{x_{r + 1}}{\koszulh{r}{\seq{x}{1}{r}}{M}} \sse
\koszulh{r}{\seq{x}{1}{r}}{M}
\]
is contained in the kernel $K$ of the canonical map (obtained from
the direct limit) $\koszulh{r}{\seq{x}{1}{r}}{M} \lra
\locoho{\f{q}_r}{r}{M}$. Recall that the cohomology module
\koszulh{r}{\seq{x}{1}{r}}{M} is the module $M/\f{q}_r M$. It
follows from the definition of the direct limit that
\[
K = \frac{\cup_{n \geq 1} \colon{(x_{1}^{n + 1}, \ldots, x_{r}^{n
+ 1})M}{M}{(x_1 \cdot \cdots \cdot x_r)^n}}{\f{q}_r M}.
\]
Furthermore,
\[
\koszulh{0}{x_{r + 1}}{\koszulh{r}{\seq{x}{1}{r}}{M}} =
\frac{\colon{\f{q}_r M}{M}{x_{r + 1}}}{\f{q}_r M}.
\]
Thus we have shown that
\[
\frac{\colon{\f{q}_r M}{M}{x_{r + 1}}}{\f{q}_r M} \sse
\frac{\cup_{n \geq 1} \colon{(x_{1}^{n + 1}, \ldots, x_{r}^{n +
1})M}{M}{(x_1 \cdot \cdots \cdot x_r)^n}}{\f{q}_r M}.
\]
Since the union on the right side is a union of increasing
submodules, there exists an integer $n$ so that
\[
\frac{\colon{\f{q}_r M}{M}{x_{r + 1}}}{\f{q}_r M} \sse
\frac{\colon{(x_{1}^{n + 1}, \ldots, x_{r}^{n + 1})M}{M}{(x_1
\cdot \cdots \cdot x_r)^n}}{\f{q}_r M}.
\]
Thus our proof is complete.
\end{proof}

\begin{proposition}\label{P:BustUpUnmixed}
Suppose that $M$ has finite local cohomologies and let \f{a} be a standard ideal for $M$.  Suppose
that $\locoho{\f{m}}{r}{M} = 0$ for some integer $r$ with $1~\leq~r~\leq~d - 1$. For any part of a
system of parameters $x_1$, \ldots, $x_r$ of $M$, if $\ideal{x}{1}{r}{A} \sse \f{a}$, then
\[
\unmixed{\ideal{x}{1}{r}{M}} = \ideal{x}{1}{r}{M} + \sum_{i =
1}^{r} \unmixed{(x_1, \ldots, \widehat{x_i}, \ldots, x_r)M}.
\]
\end{proposition}

\begin{proof}
Let $x_1$, \ldots, $x_r$ be as in the statement of the
proposition.  For each integer $i$ with $1 \leq i \leq r$, set
$\f{q}_i = (x_1, \ldots, \widehat{x_i}, \ldots, x_r)A$, where
$\f{q}_i = 0$ if $r = 1$. Set $\f{q} = \ideal{x}{1}{r}{A}$. Choose
an element $x_{r+1}$ in \f{a} so that $x_1$, \ldots, $x_{r+1}$ is
again part of a system of parameters for $M$.  According to
Lemma~\ref{L:ColonToUnmixeds}, for any integer $n > 0$ we have
\[
\colon{(x_{1}^{n + 1}, \ldots, x_{r}^{n + 1})M}{M}{(x_{1} \cdot
\cdots \cdot x_{r})^n} = \f{q}M + \sum_{i = 1}^{r}
\unmixed{\f{q}_i M}.
\]
From this and Lemma~\ref{L:UnmixedInsideColon}, it follows that
\[
\colon{\f{q}M}{M}{x_{r + 1}} \sse \f{q}M + \sum_{i = 1}^{r}
\unmixed{\f{q}_i M}.
\]
From Lemma~\ref{L:UnmixedEqualsColon}, we have that
\[
\unmixed{\f{q}M} \sse \f{q}M + \sum_{i = 1}^{r} \unmixed{\f{q}_i
M}.
\]

Using Lemma~\ref{L:UnmixedEqualsColon} again, we see that
$\unmixed{\f{q}_i M} \sse \unmixed{\f{q}M}$ for each integer $i$
with $1 \leq i \leq r$. This completes the proof.
\end{proof}

\subsection{The intersection of an unmixed component with a system of parameters}

The proof of the following proposition is inspired by
\cite[Lemma~2.2]{Kawasaki}, a similar result that is stated in
terms of $d$-sequences.

\begin{proposition}\label{P:UnmixedIntersectParameterIdeal}
Suppose $M$ has finite local cohomologies and let \f{a} be a
standard ideal for $M$.  Let $x_1$, \ldots, $x_r$, \ldots, $x_{r +
n}$ $(0 \leq r < d$, $n \geq 0)$ be part of a system of parameters
for $M$.  If $\ideal{x}{1}{r + n}{A} \sse \f{a}$ then
\[
\unmixed{\ideal{x}{1}{r}{M}} \intersect \ideal{x}{1}{r + n}{M} =
\ideal{x}{1}{r}{M}.
\]
\end{proposition}

\begin{proof}
We go by induction on $n$.  The result is trivial when $n = 0$.

Suppose $n > 0$.  Let $a$ be in $\unmixed{\ideal{x}{1}{r}{M}}
\intersect \ideal{x}{1}{r + n}{M}$.  Write $a = \sum_{i = 1}^{r +
n} x_{i} a_{i}$ with each $a_i$ in $M$. Since
\[
\unmixed{\ideal{x}{1}{r}{M}} = \colon{\ideal{x}{1}{r}{M}}{M}{x_{r
+ n}},
\]
we multiply $a$ by $x_{r + n}$ and obtain
\[
\sum_{i = 1}^{r + n} x_{i} x_{r + n} a_{i} \in \ideal{x}{1}{r}{M}.
\]
Examining the highest term, we see that $x_{r + n}^{2} a_{r + n}$ is
in the submodule \ideal{x}{1}{r + n - 1}{M}.  Since \f{a} is a
standard ideal for $M$, this implies that
\[
x_{r + n} a_{r + n} \in \ideal{x}{1}{r + n - 1}{M}.
\]
Recalling our expression for $a$, we see that $a$ is in
\ideal{x}{1}{r + n - 1}{M}.  Thus $a$ is in
$\unmixed{\ideal{x}{1}{r}{M}} \intersect \ideal{x}{1}{r + n -
1}{M}$.  By the induction hypothesis, we see that $a$ is in
\ideal{x}{1}{r}{M}, and the proof is complete.
\end{proof}

\subsection{Local cohomology as a quotient of an unmixed component}

\begin{lemma}\label{L:RegularSequence}
Suppose $M$ has finite local cohomologies and let \f{a} be a
standard ideal for $M$.  Set $t = \depth M$.  If $x_1$, \ldots,
$x_t$ is part of a system of parameters for $M$ contained in
\f{a}, then $x_1$, \ldots, $x_t$ is a regular sequence on $M$.
\end{lemma}

\begin{proof}
There is nothing to show if $t = 0$. Suppose $t
> 0$. Since $\colon{0}{M}{x_1} = \colon{0}{M}{\f{a}}$,
the submodule \colon{0}{M}{x_1} has finite length. Since $M$ has
depth $t
> 0$, it must be that $\colon{0}{M}{x_1} = 0$, because $M$ has no
nonzero submodule of finite length.  Hence $x_1$ is regular on
$M$.

If $t = 1$, we are done.  Otherwise, we go modulo the regular
element $x_1$ and continue as above.  Specifically, since
$\colon{x_1 M}{M}{x_2} = \colon{x_1 M}{M}{\f{a}}$, \f{a} kills the
module \colon{0}{M/x_1 M}{x_2}, whence this module has finite
length.  Since $M/x_1 M$ has depth $t - 1 > 0$, we see that
$\colon{0}{M/x_1 M}{x_2} = 0$, hence $x_2$ is regular on $M/x_1
M$.  We continue in this manner to complete the proof.
\end{proof}

\begin{remark}\label{R:StandardKills}
Suppose $M$ has finite local cohomologies and let \f{a} be a
standard ideal for $M$.  Let $x_1$, \ldots, $x_d$ be a system of
parameters for $M$ and let \f{q} denote the ideal they generate.
According to \cite[Theorem and Definition~17,
p.~261]{Stuckrad.Vogel}, if $\f{q} \sse \f{a}$, then we have
\[
\f{q} \locoho{\f{m}}{i}{M/\ideal{x}{1}{j}{M}} = 0
\]
for all $i$ and $j$ with $j \geq 0$ and $0 \leq i < d - j$.
\end{remark}

\begin{lemma}\label{L:GoModOneElement}
Suppose $x$ is an element of $A$ that is regular on $M$.  If $x$ annihilates \locoho{\f{m}}{r}{M}
and if $\locoho{\f{m}}{r - 1}{M} = 0$, then
\[
\locoho{\f{m}}{r}{M} \cong \locoho{\f{m}}{r - 1}{M/xM}.
\]
\end{lemma}

\begin{proof}
From the short exact sequence
\[
0 \lra M \lra M \lra M/xM \lra 0
\]
induced by multiplication by $x$ on $M$, we obtain the following
isomorphism from the long exact sequence for local cohomology:
\[
0 \lra \locoho{\f{m}}{r - 1}{M/xM} \lra \locoho{\f{m}}{r}{M} \lra
0.
\]
The zero on the left side comes from the fact that
$\locoho{\f{m}}{r - 1}{M} = 0$; the zero on the right side is due
to the fact that $x \locoho{\f{m}}{r}{M} = 0$.
\end{proof}

\begin{lemma}\label{L:GoModSeveralElements}
If $x_1$, \ldots, $x_r$ is a regular sequence on $M$ such that
$\ideal{x}{1}{r}{A}$ annihilates \locoho{\f{m}}{r}{M}, then
\[
\locoho{\f{m}}{r}{M} \cong
\locoho{\f{m}}{0}{M/\ideal{x}{1}{r}{M}}.
\]
\end{lemma}

\begin{proof}
We go by induction on $r$.  If $r = 1$ then we are done by
Lemma~\ref{L:GoModOneElement}.

Now assume $r > 1$.  By Lemma~\ref{L:GoModOneElement}, we have an
isomorphism $\locoho{\f{m}}{r}{M} \cong \locoho{\f{m}}{r - 1}{M/x_1
M}$. Set $N = M/x_1 M$.  By the induction hypothesis applied to $N$
and the sequence $x_2$,~\ldots,~$x_r$, we have $\locoho{\f{m}}{r -
1}{N} \cong \locoho{\f{m}}{0}{N/\ideal{x}{2}{r}{N}}$.  We complete
the proof using this and the previous isomorphism.
\end{proof}

\begin{proposition}\label{P:Locoho=UnmixedModParameter}
Suppose $M$ has finite local cohomologies and let \f{a} be a
standard ideal for $M$.  Let $x_1$, \ldots, $x_r$ $(0 \leq r \leq
d - 1)$ be part of a system of parameters for $M$ and suppose
$\depth M \geq r$.  If $\ideal{x}{1}{r}{A} \sse \f{a}$ then
\[
\locoho{\f{m}}{r}{M} \cong
\unmixed{\ideal{x}{1}{r}{M}}/\ideal{x}{1}{r}{M}.
\]
\end{proposition}

\begin{proof}
By Lemma~\ref{L:RegularSequence}, \seq{x}{1}{r} forms a regular
sequence on $M$.  By Remark~\ref{R:StandardKills} we have
$\ideal{x}{1}{r}{\locoho{\f{m}}{j}{M}} = 0$ for all $j$ such that
$0 \leq j < d$, so from Lemma~\ref{L:GoModSeveralElements} we
obtain the isomorphism $\locoho{\f{m}}{r}{M} \cong
\locoho{\f{m}}{0}{M/\ideal{x}{1}{r}{M}}$.  We complete the proof
with an application of Remark~\ref{R:Unmixed}.
\end{proof}

\subsection{Towards the directness of a sum}

This subsection is dedicated to the proof of the following
proposition.

\begin{proposition}\label{P:SumIsDirect}
Suppose $M$ has finite local cohomologies and let \f{a} be a
standard ideal for $M$.  Let \seq{x}{1}{d} be a system of
parameters for $M$ and suppose $\ideal{x}{1}{d}{A} \sse \f{a}$. If
$M$ has positive depth and $r$ is an integer with $1 \leq r \leq
\depth M$, then
\[
\unmixed{\ideal{x}{1}{r}{M}}
    \ \intersect\
    \sum_{i = 1}^{r}
        \unmixed{
            \left(
                x_{1}, \ldots, \leaveout{x_{i}}, \ldots, x_{d}
            \right)
        M}
    \ \sse\
    \ideal{x}{1}{r}{M}.
\]
\end{proposition}

\begin{proof}
Let
\[
N = \unmixed{\ideal{x}{1}{r}{M}}
    \ \intersect\
    \sum_{i = 1}^{r}
        \unmixed{
            \left(
                x_{1}, \ldots, \leaveout{x_{i}}, \ldots, x_{d}
            \right)
        M}.
\]
Since $\unmixed{\ideal{x}{1}{r}{M}} =
\colon{\ideal{x}{1}{r}{M}}{M}{x_{d}}$, we see that the submodule
$x_{1} \cdots x_{r} x_{d} N$ is contained in
\[
x_{1} \cdots x_{r} \ideal{x}{1}{r}{M}
    \ \intersect\
    \sum_{i = 1}^{r}
        x_{1} \cdots \leaveout{x_{i}} \cdots x_{r} x_{d}
            \left(
                x_{1}, \ldots, \leaveout{x_{i}}, \ldots, x_{d}
            \right) M.
\]
By considering $x_{d} M$ separately, we see that the previous
expression equals
\begin{align*}
&x_{1} \cdots x_{r} \ideal{x}{1}{r}{M}
    \intersect
    \left[
        \sum_{i = 1}^{r}
            x_{1} \cdots \leaveout{x_{i}} \cdots x_{r} x_{d}
            \left(
                x_{1}, \ldots, \leaveout{x_{i}}, \ldots, x_{d-1}
            \right) M \right. \\
        &\quad \quad \quad \quad \quad \quad \quad \quad \quad \quad \quad \quad\left. \ +\
        \sum_{i = 1}^{r}
            x_{1} \cdots \leaveout{x_{i}} \cdots x_{r} x_{d}^{2} M
    \right].
\end{align*}
We enlarge the product on the left and see that the last expression
is contained in
\begin{align*}
&\left(
    x_{1}^{2}, \ldots, x_{r}^{2}, \seq{x}{r + 1}{d - 1}
 \right) M\\
&\quad
    \intersect
    \left[
        x_{d} \sum_{i = 1}^{r}
            x_{1} \cdots \leaveout{x_{i}} \cdots x_{r}
                \left(
                    x_{1}, \ldots, \leaveout{x_{i}}, \ldots, x_{d - 1}
                \right) M \right. \\
        &\quad \quad \quad \left. \ +\
        x_{d}^{2} \sum_{i = 1}^{r}
            x_{1} \cdots \leaveout{x_{i}} \cdots x_{r} M
    \right].
\end{align*}
Since the first summand in the right side of the intersection is
contained in the left side of the intersection, the previous
expression equals
\begin{equation}\label{E:Sum}
\begin{aligned}
&x_{d} \sum_{i = 1}^{r}
    x_{1} \cdots \leaveout{x_{i}} \cdots x_{r}
        \left(
            x_{1}, \ldots, \leaveout{x_{i}}, \ldots, x_{d - 1}
        \right) M\\
&\qquad +
    \left(
        x_{1}^{2}, \ldots, x_{r}^{2}, \seq{x}{r + 1}{d - 1}
    \right) M
    \ \intersect\
    x_{d}^{2} \sum_{i = 1}^{r}
        x_{1} \cdots \leaveout{x_{i}} \cdots x_{r} M.
\end{aligned}
\end{equation}
Since \f{a} is a standard ideal, anything multiplying $x_{d}^{2}$ into the submodule
\[
\left( x_{1}^{2}, \ldots, x_{r}^{2}, \seq{x}{r + 1}{d - 1} \right) M
\]
also multiplies $x_{d}$ into this submodule.  (We use this sort of maneuver several more times in
the course of this proof.) Thus expression \eqref{E:Sum}, and hence $x_{1} \cdots x_{r} x_{d} N$,
is contained in
\[
x_{d} \left( x_{1}^{2}, \ldots, x_{r}^{2}, \seq{x}{r + 1}{d - 1}
\right) M.
\]
Since the depth of $M$ is positive and $x_{d}$ is part of a system
of parameters for $M$ that is contained in \f{a}, $x_{d}$ is
regular on $M$.  Hence
\[
x_{1} \cdots x_{r} N \sse \left( x_{1}^{2}, \ldots, x_{r}^{2},
\seq{x}{r + 1}{d - 1} \right) M.
\]

Now we use a similar technique that is somewhat simpler.  We have that the submodule $x_{1} \cdots
x_{r} x_{d - 1} N$ is contained in
\[
x_{1} \cdots x_{r} \ideal{x}{1}{r}{M}
    \intersect
    \left[
        x_{d - 1}
        \left(
            x_{1}^{2}, \ldots, x_{r}^{2}, x_{r + 1}, \ldots, x_{d - 2}
        \right) M
        + x_{d - 1}^{2} M
    \right]
\]
Enlarging the left side of the intersection, we obtain
\begin{align*}
&\left(
    x_{1}^{2}, \ldots, x_{r}^{2}, x_{r + 1}, \ldots, x_{d - 2}
\right) M\\
    &\quad \quad \quad \quad \intersect
    \left[
        x_{d - 1}
        \left(
            x_{1}^{2}, \ldots, x_{r}^{2}, x_{r + 1}, \ldots, x_{d - 2}
        \right) M
        + x_{d - 1}^{2} M
    \right]
\end{align*}
This equals
\begin{align*}
&x_{d - 1}
    \left(
        x_{1}^{2}, \ldots, x_{r}^{2}, x_{r + 1}, \ldots, x_{d - 2}
    \right) M\\
    &\quad \quad \quad \quad \quad \ +\ \left(
        x_{1}^{2}, \ldots, x_{r}^{2}, x_{r + 1}, \ldots, x_{d - 2}
    \right) M
    \intersect
    x_{d - 1}^{2} M.
\end{align*}
Using the fact that \f{a} is a standard ideal (as above), we see
that this last expression is just $x_{d - 1} \left( x_{1}^{2},
\ldots, x_{r}^{2}, x_{r + 1}, \ldots, x_{d - 2} \right) M$.  Since
$x_{d - 1}$ is regular on $M$, we now see that
\[
x_{1} \cdots x_{r} N \sse \left( x_{1}^{2}, \ldots, x_{r}^{2},
x_{r + 1}, \ldots, x_{d - 2} \right) M.
\]

Continuing with the same procedure, we eventually arrive at the
situation where \( x_{1} \cdots x_{r} N \sse \left( x_{1}^{2},
\ldots, x_{r}^{2} \right) M \).  Since $M$ has depth at least $r$
and \seq{x}{1}{r} is part of a system of parameters for $M$
contained in a standard ideal, \seq{x}{1}{r} is actually a regular
sequence on $M$ (in any order).  A simple proof using this fact
shows that $N \sse \ideal{x}{1}{r}{M}$, as desired. This completes
the proof.
\end{proof}

\section{Examples}

The first subsection in this section gives examples of rings to which our main
theorem applies.  The second subsection gives an example to show that eventual
constant index of reducibility does not imply finite local cohomologies in
general.

\subsection{Some Rings Having Finite Local Cohomologies}

It is easy to produce examples for the Main Theorem of this paper in the case $d=1$: Any
Noetherian local ring of dimension one will do.  More care is required in case $d>1$; the
proposition we quote below provides us with plenty of examples.

\begin{proposition}\cite[Theorem~A]{Evans.Griffith}\hfill\\
Let $k$ be an infinite field and let $R = k [ \seq{X}{1}{n} ]$ be the full ring of
polynomials over $k$ with $n$ at least $4$. Suppose that $2 \leq t_{1} < \cdots < t_{s}
\leq n-2$ is a sequence of integers and that $\seq{L}{1}{s}$ are graded $R$-modules of
finite length.  Then there is a graded prime ideal $I$ such that \locoho{\f{m}}{i}{R/I}
is zero unless $i = t_{1} - 1, \ldots, t_{s} - 1, n - 2$, while \locoho{\f{m}}{t_{j} -
1}{R/I} is isomorphic to $L_{j}$ for $j = 1, \ldots, s$.  Moreover, if $t_{1}$ is at
least $3$, then $R/I$ may also be taken to be a normal domain.
\end{proposition}

In order to obtain examples relevant to this paper, we take $s=1$
and $t_{1} = r+1$, where $r$ is the integer mentioned in the Main
Theorem. We may then choose $L_1$ to be any $R$-module of finite
length; $L_1$ will be the $r$th local cohomology module of $R/I$ and
all other local cohomology modules of $R/I$ will be zero, except of
course for the highest, at position $n-2$.  To obtain local
examples, we can localize at the maximal homogeneous ideal
\ideal{X}{1}{n}{R}.

Now we turn to a more concrete example which is suggested to us by
Goto. Let $k$ be a field and let $s$ and $t$ be indeterminates.  In
what follows, we work with a graded ring over a field instead of a
local ring; to obtain a local example, we can localize at the
maximal homogeneous ideal.  Set $A = k [ s^5, s^4t, st^4, t^5 ]$ and
let $\f{m}$ denote the maximal homogeneous ideal.  Then $A$ is a
two-dimensional graded domain, where the grading is done by total
degree.  We will show that $A$ has finite local cohomologies, that
$\f{m}$ does not kill the first local cohomology module of $A$, and
that $A$ does not have constant index of reducibility of parameter
ideals.

To see that $A$ has finite local cohomologies, according to
Proposition~\ref{P:FLC=standardideal}, it suffices to produce a standard system of
parameters.  Let $x = s^{10}$ and $y = t^{10}$; we claim that $\{ x, y \}$ is a standard
system of parameters. Let $\f{q} = (x,y)A$.  According to \cite[Theorem and
Definition~17, p.~261]{Stuckrad.Vogel}, it suffices to show that each of the four
sequences $\{ x, y \}$, $\{ x^2, y \}$, $\{ x, y^2 \}$, and $\{ x^2, y^2 \}$ is
$\f{q}$-weak.  Since $x$ and $y$ are regular elements, we only need to check
$\colon{xA}{A}{y} = \colon{xA}{A}{\f{q}}$, $\colon{x^2A}{A}{y} = \colon{x^2A}{A}{\f{q}}$,
$\colon{xA}{A}{y^2} = \colon{xA}{A}{\f{q}}$, and $\colon{x^2A}{A}{y^2} =
\colon{x^2A}{A}{\f{q}}$.  Note that these ideals are monomial ideals, so they can be
computed directly.  In particular, we have
\begin{align*}
   & \colon{xA}{A}{y^2} = \colon{s^{10}A}{A}{t^{20}} = s^{10}A + (s^{13}t^2, s^{12}t^3)A\\
   & \colon{x^2A}{A}{y^2} = \colon{s^{20}A}{A}{t^{20}}
        = s^{20}A + (s^{23}t^2, s^{22}t^3)A
\end{align*}
It is not difficult to see $s^{13}t^2 \cdot t^{10} = s^{13}t^{12} = s^{10} \cdot s^3t^{12}$ and
$s^{12}t^3 \cdot t^{10} = s^{12}t^{13} = s^{10} \cdot s^2t^{13}$ are contained in $xA = s^{10}A$,
so we have $s^{13}t^2, s^{12}t^3 \in \colon{xA}{A}{y}$.  Thus $\colon{xA}{A}{y^2} =
\colon{xA}{A}{y} = \colon{xA}{A}{\f{q}}$.  Similarly, we see that $s^{23}t^2$ and $s^{22}t^3$ are
in \colon{x^2A}{A}{y}, so $\colon{x^2A}{A}{y^2} = \colon{x^2A}{A}{y}$.  Moreover, we also see that
$s^{23}t^2$ and $s^{22}t^3$ are in \colon{x^2A}{A}{x}.  Thus we have $\colon{x^2A}{A}{y^2} =
\colon{x^2A}{A}{y} = \colon{xA}{A}{\f{q}}$.

Since $\{ x, y\}$ is a standard system of parameters, by \cite[Theorem and Definition~17,
p.~261]{Stuckrad.Vogel}, they each kill \locoho{\f{m}}{1}{A}.  Moreover, since the ideal
$xA$ is unmixed up to $\f{m}$-primary components, we have $\locoho{\f{m}}{1}{A} \cong
\unmixed{xA}/xA$.  On the other hand, by Lemma~\ref{L:UnmixedEqualsColon}, we see that
$\unmixed{xA} = \colon{xA}{A}{y} = (s^{10}, s^{12}t^3, s^{13}t^2)A$; thus
$\locoho{\f{m}}{1}{A} \cong (s^{10}, s^{12}t^3, s^{13}t^2)A / s^{10}A$.  From this we see
that \locoho{\f{m}}{1}{A} is not killed by \f{m} but it is killed by $\f{m}^2$.
Furthermore, we compute the socle of the first local cohomology by calculating
$\colon{s^{10}}{A}{\f{m}} = (s^{10}, s^{13}t^7, s^{17}t^3)A$.  Hence the socle dimension
of \locoho{\f{m}}{1}{A} is two.

Now we use the $\check{\mathrm{C}}$ech complex to see that the socle dimension of
\locoho{\f{m}}{2}{A} is four.  According to the $\check{\mathrm{C}}$ech complex, we have
\[
\locoho{\f{m}}{2}{A} \cong \frac{A \left[ 1/s^5t^5 \right]}{A[1/s^5] + A[1/t^5]}.
\]
We examine the components of this expression and find that
\begin{align*}
& A[1/s^5] = k[t/s] \left[ s^5, 1/s^5 \right] = \sum_{n \in \mathbb{Z}} k[t/s]s^{5n},\\
& A[1/t^5] = k[s/t] \left[ t^5, 1/t^5 \right] = \sum_{n \in \mathbb{Z}} k[s/t]t^{5n},
\end{align*}
and
\[
A \left[ 1/s^5t^5 \right] = \sum_{n \in \mathbb{Z}} \sum_{\substack{\alpha + \beta = 5n \\
\alpha, \beta \in \mathbb{Z}}} ks^{\alpha}t^{\beta}.
\]
Notice that $1/st^4$, $1/s^2t^3$, $1/s^3t^2$, and $1/s^4t$ are linearly independent and contained
in the socle.  Now we just need to show that nothing outside of their span is in the socle.  To see
this, suppose $\alpha + \beta \geq 10$ and $\alpha$, $\beta \geq 1$. We show that
$1/s^{\alpha}t^{\beta}$ is not in the socle.  It must be true that either $\alpha - 4 > 0$ or
$\beta - 4 > 0$.  By symmetry, we may assume without loss of generality that $\alpha - 4 > 0$. If
$\beta = 1$, then $\alpha \geq 9$, and the element $s^5 (1/s^{\alpha}t^{\beta}) = 1/s^{\alpha -
5}t$  is not in $A[1/s^5] + A[1/t^5]$, since its denominator is not a pure power of $s$ or $t$.
Thus we see that when $\beta = 1$, the element $1/s^{\alpha}t^{\beta}$ is not in the socle. If
$\beta > 1$, then $s^4t (1/s^{\alpha}t^{\beta}) = 1/s^{\alpha - 4}t^{\beta - 1}$ is not in
$A[1/s^5] + A[1/t^5]$, and we see that $1/s^{\alpha}t^{\beta}$ is not in the socle. Hence
$\socdim{\locoho{\f{m}}{2}{A}} = 4$.

According to the Main Theorem, the index of reducibility of parameter ideals for $A$
contained in high powers of the maximal ideal is $2 \times 2 + 4 = 8$.  Using the fact
that the index of reducibility of a parameter ideal \f{q} equals the socle dimension of
$A/\f{q}$, we find that $\index{A}{(s^5, t^5)A}{A} = 2$, $\index{A}{(s^{10}, t^{10})A}{A}
= 4$, and $\index{A}{(s^{15}, t^{15})A}{A} = 8$, thus $A$ does not have constant index of
reducibility.  More precisely, with direct computation, we have
\begin{align*}
    & \colon{(s^5, t^5)A}{A}{\f{m}} = (s^5, t^5, s^{12}t^3, s^3t^{12})A \\
    & \colon{(s^{10}, t^{10})A}{A}{\f{m}} = (s^{10}, t^{10}, s^{17}t^3, s^{13}t^7,
        s^7t^{13}, s^3t^{17})A \\
    & \colon{(s^{15}, t^{15})A}{A}{\f{m}} = (s^{15}, t^{15}, s^{22}t^3, s^{18}t^7,
        s^{14}t^{11}, s^{13}t^{12},\\
    &\quad \quad \quad \quad \quad \quad \quad \quad \quad \quad \quad \quad \quad s^{12}t^{13}, s^{11}t^{14}, s^7t^{18}, s^3t^{22})A,
\end{align*}
so
\begin{align*}
    & \index{A}{(s^5, t^5)A}{A} = \socdim{A/(s^5, t^5)A} = 2 \\
    & \index{A}{(s^{10}, t^{10})A}{A} = \socdim{A/(s^{10}, t^{10})A} = 4 \\
    & \index{A}{(s^{15}, t^{15})A}{A} = \socdim{A/(s^{15}, t^{15})A} = 8
\end{align*}
We can also confirm certain assertions in this example by using the computer program Macaulay 2
\cite{Macaulay2}.

\begin{remark}
The preceding example shows that parameter ideals contained in a standard ideal need not
have the same index of reducibility.  Indeed, $(s^{10}, t^{10})A$ and $(s^{15}, t^{15})A$
are both contained in the standard ideal $(s^{10}, t^{10})A$, but they have different
indexes of reducibility.
\end{remark}

\subsection{Eventual Constancy Does Not Imply Finite Local Cohomologies}

Eventual constant index of reducibility does not imply finite local
cohomologies.  In fact, even constant index of reducibility does not
imply finite local cohomologies, as one can see from
\cite[Example~4.7]{Goto.Suzuki}. In this section we provide another
example by slightly generalizing and reexamining an example from
Section 5 of \cite{Goto.SakuraiI}.

\begin{example}
Let $n \geq 2$ and let $S$ be an $n+1$-dimensional regular local ring whose
maximal ideal is $(x, y_{1}, y_{2}, \ldots, y_{n})S$.  Let $A = S/(x y_{1}, x
y_{2}, \ldots, x y_{n})S$ and let $\f{m} = (x, y_{1}, y_{2}, \ldots, y_{n})A$,
the maximal ideal of $A$. Then $A$ is a Noetherian local ring with dimension
$n$ and depth 1 such that every parameter ideal contained in $\f{m}^{2}$ has
index of reducibility 2, and $A$ does not have finite local cohomologies.
\end{example}

We note that this example does not have constant index of reducibility, since
the parameter ideal $(y_{1} - x, y_{2}, \ldots, y_{n})A$ is irreducible.

\begin{proof}
Let $\f{p} = (y_{1}, y_{2}, \ldots, y_{n})A$.  The ring $A$ does not have
finite local cohomologies since its minimal primes $xA$ and \f{p} have
different dimensions:  $\Dim A/(x) = n$ and $\Dim A/\f{p} = 1$
\cite[Proposition 16, p. 260]{Stuckrad.Vogel}.

Let \f{q} be a parameter ideal for $A$ contained in $\f{m}^{2}$.  Consider the
following exact sequence of $A$-modules:
\[
\xymatrix@1{
 0 \ar[r] &
A/\f{p} \ar[r]^{\cdot x} & A \ar[r] & A/xA \ar[r] & 0
  }.
\]
Since $A/xA \cong S/xS$, $A/xA$ is a Cohen-Macaulay $A$-module of dimension
$n$.  Since \f{q} is still a parameter ideal for $A/xA$, it is generated by an
$A/xA$-regular sequence, so $\tor{A}{1}{A/\f{q}}{A/xA} = 0$.  We tensor with
$A/\f{q}$ and use fundamental isomorphisms to obtain an exact sequence
\[
\xymatrix@1{
 0 \ar[r] &
A/(\f{p} + \f{q}) \ar[r]^-{\cdot x} & A/\f{q} \ar[r] & A/(xA + \f{q}) \ar[r] &
0
  }.
\]
Now we apply the socle functor:
\[
\xymatrix@1{
 0 \ar[r] &
\soc{A/(\f{p} + \f{q})} \ar[r]^-{\cdot x} & \soc{A/\f{q}} \ar[r] & \soc{A/(xA +
\f{q})}
  }.
\]
Since $A/\f{p}$ is a DVR and $A/(\f{p} + \f{q})$ has finite length,
it follows that $\socdim{A/(\f{p} + \f{q})} = 1$.  Since $A/xA$ is a
type 1 Cohen-Macaulay $A$-module and \f{q} is a parameter ideal for
$A/xA$, $\socdim{A/(xA + \f{q})} = 1$.  Thus, as soon as we show the
map $A/\f{q} \lra A/(xA + \f{q})$ is surjective on the socles, we
will have $\index{A}{\f{q}}{A} = 1 + 1 = 2$, and our proof will be
complete.

Write $\f{q} = (f_{1}, f_{2}, \ldots, f_{n})A$.  We will modify
these elements $f_{i}$ during the course of the proof.  Use \ol{*}
to denote reduction  modulo \f{p}.  Since \ol{A} is a DVR, we know
that one of the ideals $\ol{f_{i}} \ol{A}$ contains the others.
Without loss of generality, we assume it is $\ol{f_{1}} \ol{A}$.
Thus there are elements $u_{i} \in A$ such that $\ol{f_{i}} =
\ol{u_{i} f_{1}}$ for $2 \leq i \leq n$.  Hence $f_{i} - u_{i} f_{1}
\in \f{p}$ for each $2 \leq i \leq n$.  Since $\f{q} = (f_{1}, f_{2}
- u_{2} f_{1}, \ldots, f_{n} - u_{n} f_{1})A$, we may replace
$f_{i}$ by $f_{i} - u_{i} f_{1}$ for $2 \leq i \leq n$ to arrive at
a situation where $\f{q} = (f_{1}, f_{2}, \ldots, f_{n})A$ with each
of $f_{2}$, \ldots, $f_{n}$ in \f{p}.

Since \f{q} is not contained in \f{p}, we know that $f_{1}$ is not in \f{p}.
Since $\f{m} = xA + \f{p}$, we may write $f_{1} = \epsilon x^{m} + g$ with
$\epsilon$ a unit in $A$ and $g \in \f{p}$. Thus $\f{q} = (x^{m} +
\epsilon^{-1} g, f_{2}, \ldots, f_{n})A$.  We replace $f_{1}$ by $\epsilon^{-1}
g$ to arrive at a situation where $\f{q} = (x^{m} + f_{1}, f_{2}, \ldots,
f_{n})A$ with each $f_{i}$ in \f{p}.  Furthermore, since $\f{q} \sse \f{m}^2$,
we know that $m \geq 2$ and each $f_{i}$ is in $\f{m}^2$.

Let $a$ be an element of $A$ that maps to a socle generator for $A/(xA +
\f{q})$; thus $a \f{m} \sse xA + \f{q}$. Our goal is to show that $a \f{m} \sse
\f{q}$.  Since the image of \f{m} and \f{p} are the same in $A/xA$, we may
assume $a \in \f{p}$.  Since $xA \cap \f{p} = 0$, $ax = 0$. It remains to see
that $a y_{i} \in \f{q}$ for each $1 \leq i \leq n$.

Since $a \f{m} \sse \f{q} + xA = (f_{1}, f_{2}, \ldots, f_{n}, x)A$, we have
equations
\[
a y_{i} = u_{i1} f_{1} + u_{i2} f_{2} + \cdots + u_{in} f_{n} + v_{i} x
\]
for $1 \leq i \leq n$ and elements $u_{ij}$, $v_{i}$ in $A$.  We see that each
element $v_{i} x$ is in $xA \cap \f{p}$, and is thus zero.

We may assume that each element $u_{ij}$ is either a unit or is actually in
\f{p}.  First we show that it is not possible for one of these first
coefficients to be a unit; we go by way of contradiction.  Suppose one of the
first coefficients $u_{i1}$ is a unit; without loss of generality, suppose it
is $u_{11}$.  Then we may solve for $f_{1}$ and write
\[
f_{1} = u_{11}^{-1} a y_{1} - u_{11}^{-1} u_{12} f_{2} - \ldots - u_{11}^{-1}
u_{1n} f_{n}.
\]
Hence
\[
\f{q} = (x^m + u_{11}^{-1} a y_{1}, f_{2}, \ldots, f_{n}).
\]
Since $u_{11}^{-1} a$ maps to a socle element of $A/(xA + \f{q})$ as well, we
may replace $a$ by $u_{11}^{-1} a$.  Now we have $\f{q} = (x^m + a y_{1},
f_{2}, \ldots, f_{n})$.

Set $R = A/(x, f_{2}, \ldots, f_{n})A$.  Then $R$ is a one-dimensional
Noetherian local ring with maximal ideal $\f{n} = \f{m}R$ and a parameter ideal
$Q = a y_{1}R$ such that the socle of $R/Q$ contains the image of $a$.
Furthermore, since $R$ is the quotient of the regular local ring $A/xA$ by
elements contained in the square of the maximal ideal, the multiplicity of $R$
is greater than 1.  Set $I = \colon{Q}{R}{\f{n}}$.  Then according to
\cite[Proposition 2.3]{Goto.SakuraiI}, we have that $\f{n} I = \f{n} Q$. Since
$a \in I$, $\f{n} a \sse \f{n} I = \f{n} Q = \f{n} a y_{1}$.  Thus $\f{n} a =
\f{n} a y_{1}$, so by Nakayama's Lemma we see that $\f{n} a = 0$.  This is a
contradiction, since $Q = a y_{1} R \sse \f{n} a$, and $Q$ is a parameter ideal
in the one-dimensional ring $R$.  Thus we see that none of the first
coefficients $u_{i1}$ can be a unit.

Now, since each of the $n$ first coefficients $u_{i1}$, $1 \leq i \leq n$, is a
nonunit, we may assume they are all in \f{p}, and we may rewrite the equations
as
\[
a y_{i} = u_{i1} (x^{m} + f_{1}) + u_{i2} f_{2} + \cdots + u_{in} f_{n}.
\]
Thus we see that $a \f{m} \sse \f{q}$, as desired.
\end{proof}




\end{document}